\documentclass[12pt]{amsart}

\usepackage{amssymb}

\usepackage{mathrsfs}

\usepackage{geometry}          

\usepackage{epstopdf}
\usepackage{amsthm}
\usepackage{enumerate}

\usepackage{dsfont}

\newtheorem{theorem}{Theorem}

\newtheorem{lemma}[theorem]{Lemma}

\newcommand{\const}{\mbox{\rm const.}}						% CONSTANT
\def\qed{\unskip\nobreak\hfill\penalty50\hskip 3pt\hbox{}\nobreak
\hfill\hbox{\vrule width 4 pt height 10 pt}}

\title{Entry and Return times distribution}
\author{Nicolai Haydn}
\date{\today}
\begin{document}
\maketitle

%\tableofcontents          

%%%%%%%%%%%%%%%%%%%%%%%%%%%%%%%%%%%%%%%%%%
%%%%%%%   INTRODUCTION
%%%%%%%%%%%%%%%%%%%%%%%%%%%%%%%%%%%%%%%%%%
\section{Introduction}

The area of dynamical systems is characterised by the study of longterm behaviour
as evidenced by the quest to find conditions under which a map has
 an invariant measures or stable and unstable manifolds. If there is an invariant 
 probability measure, then sensitivity of orbits on initial data expresses itself
 in ergodicity and mixing properties. One way to measure this is to look 
 at the decay of correlation functions where two functions, one of them 
 pulled back, are jointly integrated. This typically requires sufficient regularity of the
 functions. In a purely measure theoretic setting we have Poincar\'e's recurrence
 theorem (see below) which states that every point returns to itself in finite time. 
 This result has subsequently been quantified in many different ways by
 Kac and others as we will describe below. This is the point of view
 we  want to adopt in this paper where we review
 some results on entry and also return times in dynamical
systems. In the process we will focus on the statistics which is dominated by
long term independence. To quantify very short term returns
where deterministic dependency dominates and which display different statistics
is also currently worked on by various people but not considered in this review.
Also not discussed here are zero entropy systems that can have many 
different `non-standard' limiting distributions. Except for the second and last sections
 the systems considered here have positive entropy.

In this paper we discuss six aspects of return and entry times statistics as follows.
  In the second section we look at general entry and return times results which
apply to any system with an invariant ergodic probability measure. In the third section
we shortly discuss different kinds of mixing which play important roles in 
distribution results most of it require symbolic coding of orbits via measurable 
partitions. In the fourth section we look at results on the first entry and return times.
We then proceed to review past results on the entry and
return times distributions which are dominated by the approach of Galves and Schmitt
from 1997. In the fifth section we review results on higher order returns which 
we categorise by the method that was used to obtain the distribution results.
In the sixth section we look at distributions near periodic orbits and in the seventh
section we look at recurrence times which measure the time it takes for a point
to return to its own neighbourhood.

Previous reviews on this topic were done by Abadi and Galves~\cite{AG} and 
Coelho~\cite{Coe}.

%%%%%%%%%%%%%%%%%%%%%%%%%%%%%%%%%%%%%%%%%%
%%%%%%%   RETURN TIMES AND INDUCES MAP
%%%%%%%%%%%%%%%%%%%%%%%%%%%%%%%%%%%%%%%%%%
\section{Return/entry times and the induced map}

Let $T$ be a map  on a space $\Omega$ and assume $\mu$ is a $T$-invariant 
probability measure on $\Omega$. For a measurable set $U\subset\Omega$ we define
the {\em return times function}
$$
\tau_U(x)=\min\{j\ge1: T^jx\in U\}
$$
for $x\in U$ and we put $\tau_U(x)=\infty$ if the forward orbit of $x$ never intersects $U$. 
We assume $\mu(U)>0$. 

\subsection{Poincar\'e recurrence theorem}

The earliest and most famous result on the return time is due to Poincar\'e~\cite{Poin} in 1890.

\begin{theorem} ({\it Poincar\'e recurrence theorem})  Let $T:\Omega\rightarrow\Omega$
and $\mu$ be a $T$-invariant probability measures. 
If $\mu(U)>0$, then $\tau_U(x)<\infty$ for almost every $x\in U$.
\end{theorem}

\noindent This theorem is easily proven as follows: Let $U_n=\bigcup_{j=n}^\infty T^{-j}U$ for the set of points $x\in\Omega$ that 
enter $U$ at least once after time $n$. Obviously $U_0\supset U_1\supset U_2\supset\cdots$,
and also $U_n=T^{-1}U_{n+1}$ which implies by the invariance of the measure that
$\mu(U_n)=\mu(T^{-1}U_{n+1})=\mu(U_{n+1})$ and consequently $\mu(U_0)=\mu(U_n)\;\forall\;n$.
Now $W=\bigcap_{n=1}^\infty U_n=\{x\in\Omega \;\mbox{enters $U$ infinitely often}\}$
and  $V=W\cap U=\{x\in U \;\mbox{enters $U$ infinitely often}\}$. Since $\mu(U_0)=\mu(U_n)$ 
we obtain that $\mu(W)=\mu(U_0)$ and since $U\subset U_0$ we conclude that 
$\mu(V)=\mu(U)$.

The recurrence statement does not extend to infinite measures as the
example  $Tx=x+1$ on the real line shows where no set of positive
measure is recurrent.

\subsection{Kac's theorem}
For $U\subset\Omega$, $\mu(U)>0$  Poincar\'e's recurrence theorem states that
 $\tau_U(x)<\infty$ for almost every $x\in U$, however
it doesn't tell us anything about how big $\tau_U$ is. Assuming ergodicity, a well known 
theorem by Kac from 1947 tells us what the average value is, i.e. the expected value of
 $\tau_U$ on $U$ and in particular also asserts that $\tau_U$ is integrable on $U$.

\begin{theorem} \cite{Kac} If $\mu$ is an ergodic $T$-invariant probability 
measure on $\Omega$ then for any $U\subset\Omega$
of positive measure one has
$$
\int_U\tau_U(x)\,d\mu(x)=1.
$$
\end{theorem}

\noindent  This theorem is a consequence of the Birkhoff ergodic theorem but can also be proven
in other ways as well. However 
in the invertible case it is easy to see why the result is true and involves
the construction of a  Kakutani tower.
Put $U_k=\{x\in U: \tau_U(x)=k\}$, $k=1,2,\dots$, for the level sets of the return time. 
Then $U=\dot{\bigcup}_{k=1}^\infty U_k$ is a disjoint union and
 the sets $T^jU_k$ for $j=0,1,\dots, k-1, k\in\mathbb{N}$, are pairwise disjoint. Since
 $\mu$ is ergodic, $\Omega=\bigcup_k\bigcup_{j=1}^{k-1}T^jU_k$ and as 
 $T$ is invertible $\mu(T^jU_k)=\mu(U_k)$. We thus obtain the statement of Kac's theorem: 
 $$
 1=\mu(\Omega)=\sum_{k=1}^\infty\sum_{j=0}^{k-1}\mu(T^{-1}U_k)=\sum_k k\mu(U_k)=\int_U \tau_U\,d\mu.
 $$
 
\vspace{3mm}

\noindent If the function $\tau_U$ is extended to the entire space $\Omega$ then we 
refer to it at the {\em entry times function}. It is not necessarily integrable over all of 
$\Omega$. In fact $\tau_U$ is integrable over $\Omega$ if and only if it is $\mathscr{L}^2$-integrable
over $U$.

\subsection{Entry and return times distributions}
Let $B\subset \Omega$ ($\mu(B)>0$) and consider the entry times function  $\tau_B(x)$ 
where $x\in\Omega$. For (parameter values) $t>0$ put
$$
F_B(t)=\mathbb{P}\left(\tau_B>\frac{t}{\mu(B)}\right)
=\mu\left(\left\{x\in\Omega: \tau_B(x)>\frac{t}{\mu(B)}\right\}\right)
$$
for the entry time distribution to $B$. Clearly, the entry times distribution $F_B(t)$ is locally constant 
on intervals of length $\mu(B)$
 and has jump discontinuities at values $t$ which are integer multiples of $\mu(B)$. 
 For any $s\in\mathbb{N}_0$ one has
$$
\{\tau_B>s+1\}=T^{-1}\{\tau_B>s\}\setminus T^{-1}B
$$ 
and consequently
 $$
 \mathbb{P}(\tau_B=s+1)=\mathbb{P}(\tau_B>s)-\mathbb{P}(\tau_B>s+1)\le\mu(B)
 $$
which shows that the jumps at the discontinuities are bounded by $\mu(B)$.

Correspondingly, for the return time function $\tau_B|B$, we call
 $$
 \tilde{F}_B(t)=\mathbb{P}_B\left(\tau_B>\frac{t}{\mu(B)}\right)
 $$
 the {\em return times distribution}.
 
 \vspace{3mm}
 
\noindent {\bf Example.} Let $\Omega$ be the shift space over an alphabet $\mathcal A$
(which is finite or countably infinite). That is $\Omega=\mathcal{A}^\mathbb{Z}$
and points $\vec{x}$ in $\Omega$ are of the form 
$\vec{x}=(\cdots x_{-2}x_{-1}.x_0x_2\cdots)$, where $x_j\in\mathcal{A}$ and the dot 
indicates the $0$th coordinate. The map on $\Omega$  is the left shift transformation $\sigma$
given by $(\sigma(\vec{x})_j=x_{j+1} \forall j$.
 If $B=U(x_0x_1\cdots x_{n-1})=\{\vec{y}\in\Omega: y_0\cdots y_{n-1}=x_0\cdots x_{n-1}\}$
is an $n$-cylinder then entry time $\tau_B(\vec{y})$ for $\vec{y}\in \Omega$
 measures the `time' it take to see the word $x_0x_1\cdots x_{n-1}$, that is
 $$
 \tau_B(\vec{y})=\min\{j\ge1: y_j\cdots y_{j+n-1}=x_0\cdots x_{n-1}\}.
 $$
 On the other hand if we restrict to $\vec{y}\in B$ then $\tau_B$ is the return time
 and measures the time it takes for the initial $n$-word to reappear. That is
  $$
 \tau_B(\vec{x})=\min\{j\ge1: x_j\cdots x_{j+n-1}=x_0\cdots x_{n-1}\}.
 $$
 The function $\tilde{F}_B(t)$ then measures the probability to see the first $n$-word
 again after rescaled time $t/\mu(B)$.
 
 \vspace{3mm}
 
\noindent The following result relates the limiting entry times distribution to the limiting return times 
 distribution. It turns out that a simple formula allows us to compute one from the other one.

 \begin{lemma}
 For any $B\subset \Omega$ ($\mu(B)>0$) let $F_B'$ be the right sided derivative of 
 the largest continuous piecewise linear function that lies below $F_B$ and is linear on
 the intervals $[k\mu(B),(k+1)\mu(B)], k\in\mathbb{N}$.
 Then
 $$
 -F_B'(t)=\tilde{F}(t).
 $$
 \end{lemma}
 
 \noindent {\bf Proof.} We have that $F_B(t)=\mathbb{P}(\tau_B>t/\mu(B))$, 
  $\tilde{F}_B(t)=\mathbb{P}_B(\tau_B>t/\mu(B))$ and for $j=0,1,2,\dots$ put 
  $A_j=\Omega\setminus T^{-j}B=T^{-j}B^c$ and ($j\le k$)
  $D_j^k=\bigcap_{\ell=j}^kA_\ell=\{x\in\Omega: T^\ell x\not\in B\;\forall\ell=j,\dots,k\}$.
  Then for any $s\in\mathbb{N}$
  $$
  \{x\in\Omega: \tau_B(x)=s\}=T^{-s}\cap D_1^{s-1}=D_1^{s-1}\setminus D_1^s
  $$
  and hence
  $$
  \mathbb{P}(\tau_B=s)=\mu\left(D_1^{s-1}\right)-\mu\left(D_1^{s}\right)
  =\mu\left(D_0^{s-2}\right)-\mu\left(D_0^{s-1}\right)
  $$
  by invariance of $\mu$ as $D_1^{s}=T^{-1}D_0^{s-1}$.
Since we also have
  $$
  \mu\left(D_0^{s-2}\right)=\mu\left(D_1^{s-1}\right)
  =\mu\left(B\cap D_1^{s-1}\right)+\mu\left(D_0^{s-1}\right),
  $$
  we consequently get
  $$
  \mu\left(\left\{x\in B: \tau_B(x)>s-1\right\}\right)
  =\mu\left(B\cap D_1^{s-1}\right)=\mathbb{P}(\tau_B=s).
  $$
  On  the other hand, if we put $s=\frac{t}{\mu(B)}$, then we can also write
\begin{eqnarray*}
  \mathbb{P}(\tau_B=s)&=&  \mathbb{P}(\tau_B>s-1)-  \mathbb{P}(\tau_B>s)\\
  &=&\mathbb{P}\left(\tau_B>\frac{t}{\mu(B)}-1\right)-\mathbb{P}\left(\tau_B>\frac{t}{\mu(B)}\right)\\
  &=&F_B(t-\mu(B))-F_B(t)\\
  &=&-\mu(B)F_B'(t-\mu(B)).
\end{eqnarray*}
Combining this with the previous identity for $\mathbb{P}(\tau_B=s)$ yields
  $$
  -F_B'(t-\mu(B))=\frac{\mu(\{x\in B: \tau_B>(t-\mu(B))/\mu(B)\})}{\mu(B)}
  =\mathbb{P}_B\left(\tau_B>\frac{t-\mu(B)}{\mu(B)}\right)
  =\tilde{F}_B(t-\mu(B)).
  $$
  \qed

\vspace{3mm}

\noindent Now let $B_n\subset\Omega$ ($\mu(B_n)>0$) be a sequence of subsets so that 
$\mu(B_n)\rightarrow0^+$ as $n\rightarrow\infty$. If the limit 
 $F(t)=\lim_{n\rightarrow\infty}F_{B_n}(t)$ exists for almost every $t\in(0,\infty)$ then
 we say $F$ is the {\em limiting entry times distribution}. Similarly 
 $\tilde{F}=\lim_{n\rightarrow\infty}\tilde{F}_{B_n}(t)$ is the {\em limiting return times distribution}
 if the limit exists almost surely. Note that the limiting entry times distribution $F$ is 
 Lipschitz continuous with Lipschitz constant $1$. The limiting return times distribution however
 does not have to have such regularity. If we apply the last lemma to the sets $B_n$ and take a limit
 then we obtain the following translation formula (from 2005).
  
 \begin{theorem} \cite{HLV}\label{HLV}
 Let $B_n\subset\Omega$ ($\mu(B_n)>0$) be a sequence of sets so that $\mu(B_n)\rightarrow0^+$.
 If the limits $F(t)=\lim_{n\rightarrow\infty} F_{B_n}(t)$, $ \tilde{F}(t)=\lim_{n\rightarrow\infty}  \tilde{F}_{B_n}(t)$
 exist (pointwise) then 
 $$
F(t)=1-\int_0^t  \tilde{F}(s)\,ds.
 $$
 \end{theorem}
 
 \noindent Since
  the limiting return distribution $\tilde{F}(t)$ is monotonically decreasing to zero,
  one sees that $F(t)$ is convex. We cannot a priori assume that $\int_0^\infty \tilde{F}\,ds$ is 
  equal to $1$. This is still an open question at this point.
  
  We also observe that the limiting entry times distribution and return times
distribution are the same only if they are exponential, that is $\tilde{F}=F$ if only if
$F(t)=\tilde{F}(t)=e^{-t}$. This result was proven independently in~\cite{HSV}
but follows very easily now from the translation formula.
A somewhat generalised version of the theorem was proven in~\cite{AS}.

  \vspace{3mm}
  
  \noindent 
  Lacroix has shown that if $F(t)$ is an eligible limiting distribution, that is 
it satisfies $F(0)=1$, is continuous, convex, monotonically decreasing on $(0,\infty)$ and 
$F(t)\rightarrow0^+$ as $t\rightarrow\infty$,
then for any ergodic $T$-invariant probability measure $\mu$ there exists a sequence of 
positive measure sets $B_n\subset\Omega$ so that $\mu(B_n)\rightarrow0$ and
such that $F(t)=\lim_{n\rightarrow\infty}F_{B_n}(t)$ for every $t\in(0,\infty)$.
Of course, the sets $B_n$ are typically pretty wild looking and in particular they are not topological
balls or cylinder sets (if there is a partition).

\subsection{Partitions} Let us assume that $\Omega$ has a measurable partition
$\mathcal{A}$ which we assume to be generating. Then
$\mathcal{A}^n=\bigvee_{j=0}^{n-1}T^{-j}\mathcal{A}$ is its $n$th join and for a point $x\in\Omega$
we put $A_n(x)$ for the unique $n$-cylinder in $\mathcal{A}^n$ which contains $x$.
Let us put $F^n_x(t)=\mathbb{P}(\tau_{A_n(x)}> t/\mu(A_n(x)))$ for the distribution of the 
entry times to the $n$-cylinder $A_n(x)$. Below we will list some results when the 
limiting distribution is know. First however we would like to state a very general
result of Downarowicz and Lacroix.

\begin{theorem}\cite{DL11} 
Let $\mu$ be an ergodic  invariant measure of positive entropy and $\mathcal{A}$ a finite
generating partition. Then
$$
\sup_{t\in\mathbb{R}}(e^{-t}-F^n_x(t))\to 0
$$
in $\mathscr{L}^1(\Omega)$ as $n\rightarrow\infty$. 
\end{theorem}

\noindent If the limiting distribution $F_x(t)=\lim_{n\to \infty}F^n_x(t)$ exists, then
this theorem states in particular that $F_x(t)\ge e^{-t}$ almost surely.

\subsection{The induced map}

For a subset $U\subset\Omega$, $\mu(U)>0$, let us denote by 
$\hat{T}=T^{\tau_U}:U\rightarrow U$ the {\em induced map}. 
By Poincar\'e's (or Kac's) theorem $\hat{T}$ exists almost everywhere. We also have
the {\em induced measure} $\hat\mu$ which is defined on $U$ by
$\hat\mu(A)=\frac{\mu(A)}{\mu(U)}$ for all measurable $A\subset U$.
Recall that $\hat\mu$ is $\hat{T}$-invariant and also that $\hat\mu$ is ergodic if 
$\mu$ is ergodic and vice versa 
under the condition that $\Omega\subset\bigcup_jT^jU$.

The following theorem from 2011 shows that a restricted system $(U,\hat{T},\hat\mu)$ has the same
limiting entry times distribution as the original system $(\Omega, T,\mu)$.

\begin{theorem} \cite{Hay}\label{induced.distribution}
 Let $\mu$ be ergodic, $U\subset\Omega$, $\mu(U)>0$. 
 Assume there exists a sequence of sets $B_n\subset U$, $\mu(B_n)\rightarrow0^+$, so that
 \begin{eqnarray*}
 F(t)&=&\lim_{n\rightarrow\infty}F_{B_n}(t),\hspace{10mm}
 F_{B_n}(t)=\mathbb{P}\left(\tau_{B_n}>\frac{t}{\mu(B_n)}\right)\\
\hat{F}(t)&=&\lim_{n\rightarrow\infty}\hat{F}_{B_n}(t),\hspace{10mm}
\hat{F}_{B_n}(t)=\mathbb{P}\left(\hat\tau_{B_n}>\frac{t}{\hat\mu(B_n)}\right)
 \end{eqnarray*}
 where
$$
\tau_B(x)>\min\{j\ge1:T^jx\in B\},\;\;\;\;\hat\tau_B(x)>\min\{j\ge1:\hat{T}^jx\in B\}
$$
 and $\hat{T}=T^{\tau_U}$ is the induced transformation on $U$. 
 
 Then $F(t)=\hat{F}(t)$ for all $t\in\mathbb{R}^+$.
\end{theorem}

 \noindent The same result holds for the return times distributions $\tilde{F}$ and its
 counterpart for the induced map $\tilde{\hat{F}}$. In fact this result was 
 for ergodic Radon measures 
 $\mu$ proven in~\cite{BSTV} in 2003 where the Lebesgue Density theorem was used
 and the limit was along metric balls $B_n$ that shrink to a point $x\in\Omega$.
 
 In conjunction with Theorem~\ref{HLV}, we see that
the limiting return times distribution of the restricted system $(U,\hat{T},\hat\mu)$ 
(for some positive measure $U\subset\Omega$) is 
the same as the limiting return times distribution of the entire system $(\Omega,T,\mu)$.

  %%%%%%%%%%%%%%%%%%%%%%%%%%%%%%%%%%%%%%%%%%%%%%%%%%%
%%%%%%%%%%%%%%%%%%%%%%%%%%%%%%%%%%%%%%%%%%%%%%%%%%%
%%%%%%%%%%%%%%%%%%%%%%%%%%%%%%%%%%%%%%%%%%%%%%%%%%%
%%%%%%   SECTION 2: MIXING PROPERTIES
%%%%%%%%%%%%%%%%%%%%%%%%%%%%%%%%%%%%%%%%%%%%%%%%%%%

\section{Mixing properties}
 
 In the previous sections we considered ergodic measures and obtained some general
 properties of the entry and return times distributions. If we wish to prove more 
 specific results then we have to make some assumptions on mixing properties
 of the measure.

Let  $\mathcal A$ be a (possibly countably infinite)
 measurable partition of $\Omega$ and denote by
${\mathcal A}^n=\bigvee_{j=0}^{n-1}T^{-j}{\mathcal A}$
its {\em $n$-th join} which also is a measurable partition of $\Omega$ for
every $n\geq1$. The atoms of ${\mathcal A}^n$ are called {\em $n$-cylinders}.
Let us put ${\mathcal A}^*=\bigcup_{n=1}^\infty{\mathcal A}^n$ for the collection of
all cylinders in $\Omega$ and put $|A|$ for the length of a
cylinder $A\in{\mathcal A}^*$, i.e.\ $|A|=n$ if $A\in{\mathcal A}^n$.

We shall assume that $\mathcal A$ is generating, i.e.\ that the atoms of
${\mathcal A}^\infty$ are single points in $\Omega$.

\subsection{Various kinds of mixing}
The following is a list of mixing properties arranged in decreasing order of strength.
Except for the $\beta$-mixing case, $U$ is always in the $\sigma$-algebra
 generated by ${\mathcal A}^n$ and $V$ lies in the $\sigma$-algebra 
generated by ${\mathcal A}^*$ (see also~\cite{Dou,Brad}).
The limiting behaviour is as the length of the `gap' $\Delta\rightarrow\infty$:
\begin{enumerate}
\item  {\em $\psi$-mixing}: 
$\displaystyle
\sup_n\sup_{U,V}\left|\frac{\mu(U\cap T^{-\Delta-n}V)}{\mu(U)\mu(V)}-1\right|=\psi(\Delta)\rightarrow0.
$
\item  {\em Left $\phi$-mixing}:
$\displaystyle
\sup_n\sup_{U,V}\left|\frac{\mu(U\cap T^{-\Delta-n}V)}{\mu(U)}-\mu(V)\right|=\phi(\Delta)\rightarrow0.
$
\item {\em $\beta$-mixing}:
$\displaystyle\sup_{n,m}
\sum_{B\in\mathcal{A}^n,C\in T^{-n-\Delta}\mathcal{A}^m}\left|\mu(B\cap C)-\mu(B)\mu(C)\right|
\rightarrow0
$
\item {\em Strong mixing} \cite{R1} (or $\alpha$-mixing):
$\displaystyle
\sup_n\sup_{U,V}\left|\mu(U\cap T^{-\Delta-n}V)-\mu(U)\mu(V)\right|=\alpha(\Delta)\rightarrow0$.
\item   {\em Uniform mixing}~\cite{R1,R2}:
$\displaystyle
\sup_n\sup_{U,V}\left|\frac1k\sum_{j=1}^k\mu(U\cap T^{-n-j}V)-\mu(U)\mu(V)\right|\rightarrow0
$ as $k\rightarrow\infty$.

\end{enumerate}
One can also have {\em right $\phi$-mixing} 
when
$
\sup_n\sup_{U,V}\left|\frac{\mu(U\cap T^{-\Delta-n}V)}{\mu(V)}-\mu(U)\right|\rightarrow0
$
as $\Delta\rightarrow\infty$.
The strongest mixing property is $\psi$-mixing and it implies all the other kinds of mixing.
Bernoulli measures are $\psi$-mixing and so are equilibrium states on Axiom A systems
for H\"older continuous potentials. The next strongest
mixing property is $\phi$-mixing, then comes $\beta$-mixing. The uniform 
mixing property is the weakest. 

\section{First return times} The first result in dynamics was by Hirata~\cite{Hirata1,Hirata2}
in 1993 and for higher order
return times by Pitskel~\cite{Pit} in 1991. The argument of Hirata uses the exponential entry times 
distribution using the Laplace transform which technically requires a fine 
analysis of the spectrum of the transfer operator on the shift space as one
restricts to the complement of a cylinder set.

\subsection{Galves-Schmitt method} 
In 1997 Galves and Schmitt~\cite{GS} used a practical approach that uses 
the mixing property to get the limiting distribution or entry times for $\psi$-mixing measures.

\begin{theorem}\cite{GS} Let $\mu$ be a $\psi$ mixing measure on a subshift such that $\psi$ is 
summable. Then there exist constants $C>1, c_1$ and $\alpha>0$  such that 
$$
\left|\mathbb{P}\left(\tau_A\ge\frac{t}{\mu(A)\lambda_A}\right)-e^{-t}\right|
\le c_1\mu(A)^\alpha
$$
for some numbers $\lambda_A\in(C^{-1},C)$ and for all cylinders $A$. 
\end{theorem}

\noindent The proof involved to subdivide the time interval $[0,N]$ into smaller portions
and then use the fact that
$\mathbb{P}(\tau_A>s/\mu(A))=1-s+o(s)$ for very small values of $s$ which are chosen
 to be $s=\mu(A)^\beta$ for some $\beta\in(0,1)$. The $\psi$-mixing property then 
 is used to show that $\mathbb{P}(\tau_A>N)$ is approximately equal 
 to $\mathbb{P}(\tau_A>s/\mu(A))^r$ where $r=N/(s/\mu(A))$ is the number of 
 small intervals needed. The `gaps' that are opened up in order to use the mixing 
 property are negligible compared to the size of the small intervals whose lengths are $s/\mu(A)$.
 Thus one obtains the exponential distribution in the limit when $\mu(A)\rightarrow0^+$.
 The drawback of this approach is that it cannot be used for higher order returns
 because one cannot control the earlier returns and still only use the mixing property.
 However, Abadi~\cite{Ab1,Aba,Ab2,Ab4} extended this method to $\phi$-mixing and even $\alpha$-mixing
 measures and in those cases also obtained error terms.

\begin{theorem}\cite{Ab2} (2006)
Let $\mu$ be an $\alpha$-mixing measure. Suppose that $\alpha(x) ² x^{-\kappa}$
($\kappa>(1+\sqrt5)/2$. Then there exists a function $\lambda: \bigcup_n\mathcal{A}^n\rightarrow(0,3]$
such that for all $A\in\mathcal{A}^n$
$$
\lim_{n\rightarrow\infty}\sup_{t\ge 0}\left|\mathbb{P}\left(\tau_A>\frac{t}{\lambda(A)\mu(A)}\right)-e^{-t}\right|=0
$$
Moreover the rate of convergence of the above limit is bounded from above by 
$e_h(A)=C_0\inf_{n\le g\le f\le1/\mu(A)}\left(f\mu(A)+\frac{g\mu(A)+\alpha(g)}{f\mu(A)}s\right),
$
 where $s =\alpha^{-1}(\mu(A)+n)$ and $C_0$ is a positive constant.
\end{theorem}

\noindent A similar result for $\phi$-mixing measures was obtained in~\cite{AV2}
and for rational maps with critical points and equilibrium states that have a `supremum gap'
it was done in~\cite{H1}.

For more specific systems a number of results were shown in~\cite{BSTV} for the first return times.
These are one-dimensional real and complex systems that have exponential  decay of correlations.
 Here the return times are for metric balls which in the case of 
 the complex maps requires a regularising approximation of the characteristic function
 for the balls by unions of cylinder sets. For interval maps that have critical points of 
 parabolic fixed points the exponential limiting distribution is shown for an induced map
 which is uniformly expanding and then by Theorem~\ref{induced.distribution} translated
 to the original map.
 
 As it appears that all known systems with positive entropy have exponentially
 distributed limiting entry/return times, Kupsa~\cite{Kup12} has recently constructed an 
 example of a positive entropy system on a three element subshift whose
 limiting entry/return times are not exponentially distributed in the limit.
 Also, Downarowics, Lacroix and Leandri~\cite{DLL} have constructed examples where the limiting 
 entry time exists along subsequences of upper density one and can be arbitrarily 
 slowly decaying to zero. Interestingly enough, it is still an open question whether
 the limiting entry times distribution $F(t)$ decays to zero if the limit $F(t)=\lim_nF_n(t)$
 exists, i.e.\ whether the sequence $F_n$ is tight. In the mentioned examples, the limit
  does not seem to exist.

  Recently Freitas, Freitas and Todd have established a relation between extremal
  values laws and entry times distribution first for absolutely continuous measures
  in~\cite{FFT1} and then for general invariant measures in~\cite{FFT2}.

 \section{Higher order return times}
 
Let $A\subset\Omega$ and $t>0$ a parameter, then we put 
$\zeta_A^t=\sum_{j=1}^{[t/\mu(A)]}\chi_A\circ T^j$ for the counting function that
counts the number of times the orbit of a point $x$ visits the set $A$ on the 
orbit segment of length $\frac{t}{\mu(A)}$. Clearly $\zeta_A^t(x)=0$ exactly if
$\tau_A(x)>\frac{t}{\mu(A)}$. If the return times $\tau_A\circ \hat{T}$ are independent 
($\hat T$ is the induced map on $A$) and are exponentially 
distributed then $\zeta_A^t$ will be Poisson distributed. However in a deterministic 
system we only get independence of entries only in the limit as their
separation goes to infinity.
 
 The first result on higher returns was due to Doeblin in 1940 for the Gauss map at the 
origin. Then there was a long gap and nothing much seems to have happened
until 1991 when several people simultaneously began to work in this area
with different methods. 
To recall Doeblin's result on the Gauss map let  $\Omega=(0,1]$ be the 
unit interval. The Gauss map is then given by
$Tx=\frac1x\;\mbox{mod}\; 1$ and is related to the continued fraction expansion
of real numbers. If $[a_0,a_1,a_2,\dots]$ is the  continued fraction expansion of
 a point $x\in\Omega$ then
$$
x=\cfrac1{a_0+\cfrac1{a_1+\cfrac1{a_2+\cdots}}}
$$
where the integers $a_j\in\mathds{N}$ are uniquely determined by $x$
and are given by $a_j=\frac1x-Tx$. The Gauss measure $\mu$ on $(0,1]$ is the 
unique absolutely continuous $T$-invariant probability measure. Its density is
 $\frac1{\log2}\frac1{1+x}$ and Doeblin~\cite{Doe} showed that for every $\theta>0$:
 $$
 \mathbb{P}\left(\left|\left\{j:1\le j\le n, a_j(x)\ge\theta n\right\}\right|=p\right)
 \to e^{-1/(\theta\log2)}\frac1{(\theta\log2)^pp!}
 $$
 as $n\to\infty$.
 Since $x=\frac1{a_j+Tx}=\frac1{a_j+\mathcal{O}(1)}$
 (as $0<Tx\le1$) we see that  a point $x\in(0,1]$ for which 
 $\left|\left\{j:1\le j\le n, a_j(x)\ge\theta n\right\}\right|=p$ visits the interval 
 $\left(0,\frac1{n\theta}\right)$ typically exactly $p$ times on the orbit segment of length $n$.
 Since $\mu\left(\left(0,\frac1{n\theta}\right)\right)=\frac1{\log2}\log(1+\frac1{\theta n})
 \approx\frac1{\log2}\frac1{\theta n}$ we can put $A_m=(0,\frac1m)$ and see that Doeblin's
 statement translates into 
 $$
 \mathbb{P}\left(\zeta_{A_m}^t=p\right) \to e^{-t}\frac{t^p}{p!}
 $$
 as $m\to\infty$. In other words, the limiting distribution of return times at the origin is
 Poissonian.

\vspace{3mm}

 There are several ways in which limiting results on higher order returns have
 been achieved. Here we will mention results that use the moment method, 
 the Chen-Stein method, direct estimates of the total variation norm and 
 a combinatorial approach.
 
 \subsection{Moment method} 
  The first dynamical result dates to 1991 and is due to Pitskel~\cite{Pit}. It uses the 
  method of moments which was also used by Wang, Tang and Wang~\cite{WTW}
  and Denker~\cite{Den} for Gibbs measures on shift spaces.
  For a set $A\subset\Omega$ and a parameter $t>0$ let us put 
$$
\zeta_A^t=\sum_{j=1}^{N}\chi_A\circ T^j
$$
for the counting function on the orbit segment of length $N=\left[\frac{t}{\mu(A)}\right]$.
The value $\zeta_A^t(x)$ counts the number of times the orbit of $x\in\Omega$ enters
$A$ up to time $N$. Now let 
$$
G_r=\{(v_1,v_2,\dots,v_r): 1\le v_1<v_2<\cdots<v_r\le N\}
$$
be the set of possible $r$-fold return times of a point to a given set within the time interval $[1,N]$.
The moment generating function for $\zeta$ is then given by
  $$
  \varphi(z)=\sum_{k=0}^\infty z^k\mathbb{P}(\zeta=k).
  $$
 A computation verifies that  
 $$
 \varphi(z)=\sum_{r=0}^\infty\frac{(z-1)^r}{r!}\mu(\xi_r),
 $$
  where $\xi_r=\zeta(\zeta-1)(\zeta-2)\cdots(\zeta-r+1)$ is the $r$-factorial moment of $\zeta$.
  It then follows that
  $$
  \frac1{r!}\mu(\xi_r)=\sum_{\vec{v}\in G_r}\mu\left(\bigcap_{j=1}^r T^{-v_j}A\right).
  $$
  If $\zeta$ is Poisson distributed with parameter $t>0$, then $\mathbb{P}(\zeta=k)=e^{-t}\frac{t^k}{k!}$,
   the generating function is $e^{t(z-1)}$ and one has $\mu(\xi_r)=t^r$.
  The advantage of this approach lies in the fact that on the right hand side of this identity 
  one can use the mixing property, as it involves the intersection of pull backs of $A$.
  One opens up `gaps' and uses rough estimates on short returns and shows that their cardinality is `small'.
  Most of the `return time patterns' $\vec{v}$ are however long and there one can use the mixing property.
  Because one has to look at arbitrarily high orders $r$ of mixing, this approach limits itself
  to $\psi$-mixing measures like Axiom A maps.
  
  Now according to a theorem of Sevast'yanov~\cite{Sev} if the moments $\mu(\xi_r)$ for $A=A_n(x)$
   converge to  $t^r$ as 
  $n\rightarrow\infty$ then $\zeta_{A_n(x)}^t$ converges in distribution to a Poisson distribution
  with parameter $t$.
In this way Pitskel proved the following result in 1991:

\begin{theorem}\cite{Pit}
Let $\mu$ be an equilibrium state for a H\"older continuous potential on a subshift
of finite type. Then for all $t>0$
$$
\mathbb{P}\left(\zeta_{A_n(x)}^t=r\right)\rightarrow\frac{t^r}{r!}e^{-t}
$$
as $n\rightarrow\infty$ almost surely in $x\in\Omega$.
\end{theorem}

\noindent He then used an approximation argument  to extend this result to metric balls
for two-dimensional toral automorphisms.

This approach was also used in~\cite{H2} to prove that the limiting return times distributions are
Poissonian for equilibrium states of rational maps. There an approximation 
argument was used to show that the return times to metric balls are in the limit Poisson distributed.
 In 2003 we went further and obtained 
that for a somewhat more general class of $\psi$-mixing maps one has Poisson distributed
entry and return times in the limit~\cite{HV1}. An improved moment method also gave error terms 
which depend on the 
rate at which $\psi(k)$ decreases to zero as $k$ goes to infinity. The key to this is to 
obtain a Sevast'yanov type theorem that allows error estimates and to show that
 approximations of the moments
translate into an approximation of the distribution to the Poisson distribution.

\subsection{Laplace transform} This is the method that was first used by Hirata~\cite{Hirata1} 
in 1993
for equilibrium states on subshifts of finite type. He proved that the first entry time
is exponentially distributed in the limit and then argued that the weak mixing
property implies that higher order returns are in the limit Poisson distributed. This requires a careful 
analysis of the transfer operator acting on the complement of the cylinder target set 
and delicate estimates on the dominant eigenvalue $e^{P(f)}$ where $P(f)$ is the 
pressure of the potential function $f$ (see e.g.~\cite{Bow}).
This approach was
more conclusively carried out by Coelho and Collet~\cite{CC00} and also in~\cite{CCC}
for measures on subshifts that have strong mixing properties.
The combinatorics involved tend to make such an approach difficult and would favour
the moment method.

\subsection{Chen-Stein method}

Let $\nu$ be a probability measure on $\mathbb{N}_0$  and denote by
$\nu_0$ the Poisson-distribution measure with mean $t$, i.e.\
$\mathbb{P}_{\nu_0}(\{k\})=\frac{e^{-t}t^k}{k!}$  $\forall k\in \mathbb{N}_0$.
If we put $\mathcal{F}$ for the set of all real-valued functions on  $\mathbb{N}_0$,
 then the Stein 
operator $\mathcal{S}:\mathcal{F}\rightarrow \mathcal{F}$ is defined by
$$
\mathcal{S}f(k)=tf(k+1)-kf(k),\quad\text{ }  \forall k\in \mathbb{N}_0,
$$
where $f\in\mathcal{F}$. The Stein equation, given by
$$
\mathcal{S}f=h-\int_{\mathbb{N}_0}h\,d\nu_0,
$$
then has for every $\nu_0$-integrable $h\in\mathcal{F}$ a solution $f$ which is unique except for
$f(0)$ which can be chosen arbitrarily~\cite{BC1}\footnote{
In fact $f$ can be computed recursively:
$$
f(k)=\frac{(k-1)!}{t^k}\sum_{i=0}^{k-1} \left( h(i)-\nu_0(h)\right)\frac{t^i}{i!}
=-\frac{(k-1)!}{t^k}\sum_{i=k}^{\infty} \left( h(i)-\nu_0(h)\right)\frac{t^i}{i!} , 
\quad\text{} \forall k\in\mathbb{N} 
$$
}.
A probability measure $\nu$ on $\mathbb{N}_0$ is Poisson (with parameter $t$)
 if and only if $\int_{\mathbb{N}_0}\mathcal{S}f\,d\nu=0$ for all bounded functions $f\in\mathcal{F}$.

If $E\subset\mathbb{N}_0$ then one can take $h=\chi_E$ the characteristic function and
estimate the distance of a probability measure $\nu$ on $\mathbb{N}_0$ from 
the  Poisson distribution $\nu_0$ as follows:
\begin{equation}\label{finalformstein}
|\nu(E)-\nu_0(E)|=\left | \int_{\mathbb{N}_0}\mathcal{S}f\,d\mu\right |
=\left | \int_{\mathbb{N}_0}\left (t f(k+1)-kf(k)\right)d\mu\right |
=\left|\mathbb{P}(\zeta\in E)-\mu_0(E)\right|
\end{equation}
where $f$ is a solution to the Stein equation. An estimate on the function $f$ then 
allows us to bound the distance between $\nu$ and $\nu_0$ using mixing properties of up
to second order. 

\subsubsection{Toral automorphisms}
In 2004 Denker, Gordin and Sharova~\cite{DGS} applied this method to 
the Haar measure of hyperbolic toral automorphism to show that 
the limiting return times distribution is Poissonian in the limit if one considers successive returns
ball-like sets $G_n$ which contract to a non-periodic point. The return sets $G_n$
lie inside balls and have a measure that is comparable. Interestingly enough, this result 
shows a dichotomy that at non-periodic points one always gets the limit Poisson
distributed return times.

\subsubsection{Mixing measures}
Abadi used in his thesis in 2001 and in~\cite{Ab3} a theorem of Arratia, Goldstein and Gordon~\cite{AGG} 
to get the Poisson limiting distribution for $\phi$-mixing systems.  Exploiting the fact that
with the Chen-Stein method only two-fold mixing is required in order to get the Poisson 
distribution, it was  shown that
$$
\mathbb{P}(\zeta^t_{A_n}(x)=k)=e^{-t}\frac{t^k}{k!}+\mathcal{O}(f(t)\epsilon(n))
$$
for an error term function $f$ which depends on the parameter $t$ and $\epsilon$ depending
on the length $n$ of the ``target string'' and the point $x\in\Omega$.
 Different error term functions $\epsilon$ are given for $\psi$-mixing
measures and $\beta$-mixing measures. In 2008 these results were improved by 
Abadi and Vergne~\cite{AV1} to $\phi$-mixing and $\alpha$-mixing systems respectively
using a different approach. Their improved result uses a more elementary method 
combined with the exponential entry time results of Abadi's earlier papers.
All thse results hold if restricted to returns only in which case one uses the induced measure on the
initial set $A_n(\cdot)$ although the error terms are slightly different.

This extends a previous result of Hirata, Saussol and Vaienti~\cite{HSV} where
error estimates for $\mathbb{P}(\zeta^t_U)-e^{-t}\frac{t^r}{r!}$ are given under
the assumption that the partition $\{U,U^c\}$ is $\phi$-mixing or $\alpha$-mixing
(different error terms). This result was then used to obtain the Poisson distribution
and rates of convergence for higher order limiting return times for 
parabolic interval maps.

Using the Chen-Stein method, Psiloyenis proved in 2008 a similar result which 
allows the return set to be a possibly infinite unions of cylinder sets. The alphabet
in this case can be countably infinite and the entropy is not required to be finite.

\begin{theorem}\cite{Psi,HP}
Let  $\mu$ be a $\phi$-mixing $T$-invariant probability measure with an at most countably
infinite partition $\mathcal{A}$.
Let $A\in\sigma(\mathcal{A}^n)$ be a finite or infinite union of $n$-cylinders
such that $|\log\mu(A)|=\mathcal{O}(n^\eta)$ and $n^\eta\phi(n)\rightarrow0$ as $n\rightarrow\infty$
for some $\eta\ge1$ and $r_A=\mathcal{O}(n)$.
Then for all $t>0$:

\noindent {\textbf{(i) Exponential mixing rate:}}
($\phi(k)=\mathcal{O}(\vartheta^{k})$, with $\vartheta\in(0,1)$) there exists $\gamma=\gamma(\vartheta)>0$
$$
\mathbb{P}\left(\zeta_A^t=r\right)=e^{-t}\frac{t^r}{r!}
+\mathcal{O}(t (t\vee1))e^{-\gamma n}.
$$

\noindent {\textbf{ (ii) Polynomial mixing rate:}}
($\phi(k)=\mathcal{O}(k^{-\beta})$ with $\beta>1+\eta$),
$$
\mathbb{P}\left(\zeta_A^t=r\right)=e^{-t}\frac{t^r}{r!}
+\mathcal{O}(t (t\vee1))n^{-(\beta-1-\eta)}.
$$
\end{theorem}

\noindent Here $r_A$ denotes the first return of $A$ to itself, i.e. $r_A=\inf_{x\in A}\tau_A(x)$ and 
the lower bound is required to control short returns and to exclude `periodic like' behaviour.
The lower cutoff for $r_A$ is not optimal but is a typical requirement to exclude `periodic looking'
target sets $A$.

Recently Kifer~\cite{Kif12} used a similar method to obtain results on the limiting distribution of
simultaneous returns. Following the approach taken by 
Abadi using a theorem of Arratia, Goldstein and Gordon~\cite{AGG}, Kifer showed that
simultaneous returns are in the limit Poissonian. For $\ell$ simultaneous returns
one lets $q_1(j)<q_2(j)<\cdots<q_\ell(j)$ be the return times so that
$g(j)=\min_k(q_k(j)-q_{k-1}(j))$ goes to infinity as $j\to \infty$. Then one forms
the counting function for simultaneous hits at times $q_k$ by putting
$$
\xi_A^t=\sum_{j=1}^{[t/\mu(A)^\ell]}\prod_{k=1}^\ell\chi_{A}\circ T^{q_k(j)}
$$
for $x\in\Omega$. For $\ell=1$ and $q_1(j)=j$ this reduces to the standard case considered
above.

\begin{theorem} \cite{Kif12}
If $\mu$ is $\psi$-mixing, then for $A\in\mathcal{A}^n$:
$$
\mathbb{P}\left(\xi_A^t=r\right)=e^{-t}\frac{t^r}{r!}+\mathcal{E}(A).
$$
The error  term $\mathcal{E}(A)$ is up to polynomial terms in $\ell, n$ (and depending on $g$)
and exponential terms in $t$ equal to $\mu(\tilde{A})+\psi(n)$, where
$\tilde{A}=A_{r'_A}(A)$ is the $r'_A$-cylinder containing $A$
where $r'_A=\min(\inf\{j\ge0: T^jA\cap A\not=\emptyset\},n)$.
\end{theorem}

\noindent In~\cite{Kif11} this result was proven without error terms.

\subsubsection{Markov towers} L-S~Young's construction of Markov towers has proven
to be a very powerful tool to obtain results on the statistical properties of maps
and the method can also be used to obtain results on the return times distribution.
Let us recall the construction from~\cite{Y2,Y3}.

For a differentiable map $F$ on a manifold $M$ one lets $\Omega_0$ be a subset of $M$
which is partitioned into sets $\Omega_{0,i}, i=1,2,\dots$ so that
there is a return time function $R:\Omega_0\rightarrow\mathbb{N}$ which is constant on the partition elements
$\Omega_{0,i}$ and which satisfies that $F^R$ maps $\Omega_{0,i}$ bijectively to the entire set $\Omega_0$.
If we put $\Omega_{j,i}=\{(x,j): x\in\Omega_{0,i}\}$ for $j=0,1,\dots, R(\Omega_{0,i})-1$
then $\Omega =\bigcup_{i=1}^\infty\bigcup_{j=0}^{R(\Omega_{0,i})-1}\Omega_{j,i}$ is called a
{\em Markov tower} for the map $T$ given by
$$
T(x,j)=\left\{\begin{array}{ll}(x,j+1)&\mbox{if $j<R(x)-1$}\\(F^R(x),0)&\mbox{if $j=R(x)-1$}\end{array}\right..
$$ 
It has the (typically) countably infinite partition $\mathcal{A}=\{\Omega_{j,i}:\; i,j\}$.

The separation function $s(x,y)$ is the smallest positive $n$ so that $(T^R)^nx$ and  $(T^R)^ny$ lie in distinct
sub-partition elements $\Omega_{0,i}$ of $\Omega_0$. Two points $x$ and $y$ in $\Omega$ 
belong to the same $N$-cylinder if and only if
 they remain together (in the same partition element) for at least $n$ iterations of $T^R$, i.e. if $s(x,y)\ge n$,
 where $N=\sum_{j=0}^{n-1}R(T^R)^j$.
 
 The space of H\"older continuous functions $\mathcal{C}_\gamma$ consists of all functions $\varphi$ 
 on $\Omega$ for which $|\varphi(x)-\varphi(y)|\le C_\varphi \gamma^{s(x,y)}$. The norm on $\mathcal{C}_\gamma$
 is $\|\varphi\|_\gamma=|\varphi|_\infty+C_\varphi$, where $C_\varphi$ is smallest possible.
  
 Let $\nu$ be a finite given `reference' measure on $\Omega$ and assume that the Jacobian
 $JT^R$ with respect to the measure $\nu$ is  H\"older continuous, that is, there exists a 
 $\gamma\in(0,1)$ so that 
 $$
 \left|\frac{JT^Rx}{JT^Ry}-1\right|\le\mbox{\rm const} \gamma^{s(T^Rx,T^Ry)}
 $$
 for all $x,y\in \Omega_{0,i}$, $i=1,2,\dots$.

 If the return time $R$ is integrable with 
 respect to  $m$ then by~\cite{Y3} Theorem~1 there exists a $T$-invariant probability measure
 $\mu$ (SRB measure) on $\Omega$ which is absolutely continuous with respect to $\nu$. Moreover the 
 density function $h=\frac{d\mu}{d\nu}=\lim_{n\rightarrow\infty}\mathcal{L}^n\lambda$ is H\"older continuous,
 where $\lambda$ can be any initial density distribution in $\mathcal{C}_\gamma$.
 The transfer operator $\mathcal{L}:\mathcal{C}_\gamma\rightarrow\mathcal{C}_\gamma$ is defined by
 $\mathcal{L}\varphi(x)=\sum_{x'\in T^{-1}x}\frac{\varphi(x')}{JT(x')}$, $\varphi\in\mathcal{C}_\gamma$,
 and has the property that $\nu$ is a fix point of its adjoint, i.e.\ $\mathcal{L}^*\nu=\nu$.
  In \cite{Y3} Theorem~2(II) the $L^1$-convergence was proven:
$$
\|\mathcal{L}^k\lambda-h\|_{L^1}\le p(k) \|\lambda\|_\gamma
$$
where the `decay function' $p(k)=\mathcal{O}(k^{-\beta})$ if the tail decays polynomially with power $\beta$,
that is if $\nu(R> j)\le \const j^{-\beta}$.
If the return times decay exponentially, i.e.\ if $\nu(R>j)\le\const \vartheta^j$ for some $\vartheta\in(0,1)$, then
there is a $\tilde\vartheta\in(0,1)$ so that
$p(k)\le \const\tilde\vartheta^k$.

\begin{theorem}\cite{HP}
As described above let $T$ be a map on the Markov Tower structure  $\Omega$
with a reference measure $\nu$ and a return time function $R$. Let  $\mu$ be the absolutely continuous 
invariant measure. Let $K$ be a constant. Then for every $A_n\in\sigma(\mathcal{A}^n)$ 
for which  $|\log\mu(A_n)|\le Kn$, $r_{A_n}>\frac{n}{2}$ and $\sup_{A_n}R\le\frac{n}2$, 
the following result holds true:

If $\nu(R>n)={\mathcal{O}}(\vartheta^n)$  ($\vartheta\in(0,1)$) or if
 $\nu(R>n)={\mathcal{O}}(n^{-\beta})$ for some $\beta>2$, then there exists $\gamma>0$ 
such that
$$
\mathbb{P}\left(\zeta_{A_n}^t=r\right)=e^{-t}\frac{t^r}{r!}
+\mathcal{O}((t\vee1)e^{-\gamma n})\quad\forall t>0 \text{ and } \forall n\in\mathbb N.
$$
\end{theorem}

\subsection{Total variation estimates} In 2010 Chazottes and Collet~\cite{CC} proved the 
Poisson distribution for Young's Markov towers in the codimension one case when
the tails decay at an exponential rate, i.e. $\nu(R>n)\sim \vartheta^n$ for a $\vartheta\in(0,1)$. 
The estimates use the decay of correlation and require that the characteristic functions
 of the metric balls be approximated by Lipschitz continuous functions. This introduces
 additional difficulties as it is necessary to control the contribution made by an 
 annulus surrounding the metric balls on which the approximating function
 interpolates between the values $0$ and $1$.

\begin{theorem}\cite{CC} Let $(M,T,\mu)$ be a non-uniformly hyperbolic dynamical system 
modelled by a Young tower whose return-time function has an exponential tail. Assume 
that the local unstable manifolds have dimension one. Denote by $\mu$ its SRB measure. 

Then there exist constants $a,b > 0$ such that for all $\rho\in (0,1)$:\\
(i) There exists a set $M_\rho$ such that
$\mu( M_\rho) =\mathcal{O}(\rho^b)$;\\
(ii) For all $x\not\in M_\rho$ one has  ($B_\rho(x)$ is the metric ball with centre $x$ and radius $\rho$)
$$
\mathbb{P}\left(\sum_{j=0}^N\chi_{B_\rho(x)}=r\right)=e^{-t}\frac{t^r}{r!}+\mathcal{O}(\rho^a)
$$
 for all $\rho\ge 0$ small enough and for every $t > 0$.
 \end{theorem}

\noindent A similar result for polynomially decaying correlations has now been proven 
by Wasilewska~\cite{Was} with error terms which are polynomially decaying in $|\log \rho|$.
 
 \section{Periodic orbits}

%%%%%%%%%%%%%%%%%%%%%%%%%%%%%%%%%%%%%%%%%%%%%
%%%%%%%%%%%%%%%%% SUBSECTION: PHI MIXING MEASURES
%%%%%%%%%%%%%%%%%%%%%%%%%%%%%%%%%%%%%%%%%%%%%

Hirata~\cite{Hirata1} and Abadi~\cite{Aba} have pointed out that at a periodic point the limiting return time is
not exponentially distributed like $e^{-t}$ but that it is a combination of a Dirac pointmass
at the origin with a rescaled exponential decay. This can easily be seen for a Bernoulli 
measure on a full two element shift $\Sigma$. On $\Sigma$ one has the left shift
$\sigma$. If $\mu$ is the Bernoulli measure for the probabilities $p, 1-p$ for some $p\in(0,1)$, 
and $\vec{x}=0^\infty$ the fixed point whose coordinates are all $0$s, then
$F_{\vec{x}}^n(t)=\mathbb{P}_{A_n(\vec{x})}(\tau_{A_n(\vec{x})}\ge t/\mu(A_n(\vec{x})))$ has the limiting
distribution $F_{\vec{x}}(t)=(1-p)e^{t/(1-p)}$ for $t>0$. Equivalently one obtains 
that the limiting distribution of $\mathbb{P}_{A_n(\vec{x})}(\tau_{A_n(\vec{x})}\ge t/[(1-p)\mu(A_n(\vec{x}))])$
is $(1-p)e^{-t}$. This is the way it is formulated in~\cite{Aba}. For an equilibrium states for a H\"older 
continuous potential $\varphi$ on a subshift of finite type this was formulated by by Hirate~\cite{Hirata1}
for arbitrary periodic points where it was shown that at a periodic point $x$ with minimal 
period $m$ 
$$
\lim_{n\to\infty}\mathbb{P}_{A_n(x)}\left(\tau_{A_n(x)}\ge \frac{t}{(1-p)\mu(A_n(x))}\right)
=(1-p)e^{-t},
$$
where $p=e^{\sum_{j=0}^{m-1}\varphi\circ T^j}$.

If $\mu$ is a $\psi$-mixing measure
then for a periodic point $x$ with minimal period $m$ the limit
$
p=\lim_{\ell\rightarrow\infty}\left|\frac1\ell \log\mu(A_{\ell m}(x))\right|
$
exists. In particular $p$ is always strictly less than $1$.
In the following we shall assume the stronger property that
$p=\lim_{n\rightarrow\infty}\frac{\mu(A_{n+m}(x))}{\mu(A_n(x))}$. This of course implies the limit
in the lemma, but we are not sure whether the reverse implication is generally true. Also put
$q_n=\sup_{\ell\ge n}\left|\frac{\mu(A_{\ell+m}(x))}{\mu(A_\ell(x))}-p\right|$. For $t>0$ and
integers $n$ we put $\zeta_n^t$ for the classical counting function on cylinder sets
$\sum_{j=0}^{N_n} \chi_{A_n(x))}\circ T^j$ with the rescaled observation time
$$
N_n=\left[\frac{t}{(1-p)\mu(A_n(x))}\right].
$$
The limiting distribution of $\zeta_n^t$ is the P\'olya-Aeppli distribution according to which the value 
$r\in\mathbb{N}_0$ is assumed with probability $e^{-t}P_r(t,p)$ where 
$$
P_r(t,p)=\sum_{j=1}^rp^{r-j}(1-p)^j\frac{t^j}{j!} \left(\begin{array}{c}r-1\\j-1\end{array}\right).
$$
To be more precise we have the following result

%%%%%%%%%%%%%%%%%%%%%%%%%%%%%%%%%%%%%%%%%%%%%
%%%%%%%%%%%%%%%%% THEOREM: PSI MIXING MEASURES
%%%%%%%%%%%%%%%%%%%%%%%%%%%%%%%%%%%%%%%%%%%%%
\begin{theorem}\label{PolyaAeppli} \cite{HV2}
Let $\mu$ be a $\psi$-mixing measure with partition $\mathcal A$
(finite or infinite), $x$ a periodic point with minimal period $m$ and $p$ and $q_n$ as above.
Then 
$$
\mathbb{P}(\zeta_n^t=r)=e^{-t}P_r+\mathcal{E}(A_n(x)).
$$
The error term $\mathcal{E}(A_n(x))$ is up to exponential terms in $t$ and rapidly decaying
terms in $r$ roughly of the form $\inf_\delta(\delta\mu(A_n(x))+\psi(\delta))+p^\frac{n}m+q_n$.
\end{theorem}

\noindent In the case of an infinite partition $\mathcal A$ no finiteness of entropy
is required.

If $\mu$ is an equilibrium state for a H\"older continuous function $f$ on an 
Axiom A space (shift space) which has the finite, generating partition $\mathcal A$ (see~\cite{Bow})
then the error term $\mathcal{E}$ can be optimised to yield, again 
up to exponential terms in $t$ and rapidly decaying terms in $r$, roughly 
$n\mu(A_n)+p^\frac{n}m$.
In~\cite{FFT3} the compound Poisson distribution for extremal values distribution and 
by extension also for the return times distribution was proven at repelling fixed points
for some non-uniformly hyperbolic systems.

 Recently in 2012 Kifer proved a general result on points where the limiting 
  distribution is not Poissonian. For $\ell$ simultaneous returns
at times $q_k(j)=d_kj$, $k=1,\dots,\ell$ where $1\le d_1<d_2<\cdots< d_\ell$ are integers
 one puts as before
$$
\xi_A^t=\sum_{j=1}^{[t/\mu(A)^\ell]}\prod_{k=1}^\ell\chi_{A}\circ T^{d_kj}.
$$
Let $r_A=\inf\{j\ge0: T^jA\cap A\not=\emptyset\}$, put
$\kappa_A=\mbox{\rm lcm}\left(r_A/\mbox{\rm gcd}(r_A,d_i): i\right)$
and define ($r=r_A$)
$$
\rho_A=\prod_{k=1}^\ell\mu_A
\left(\bigcap_{i=0}^{[(n+d_k\kappa)/r]}T^{-ir}A\cap T^{-[n/r]}A_{n-r[n/r]}(A)\right),
$$
where $A_{n-r[n/r]}(A)$ denotes the unique $n-r[n/r]$-cylinder containing $A$.

\begin{theorem} \cite{Kif12}
If $\mu$ is $\psi$-mixing, then for every $A\in\mathcal{A}^n$ and $n$ which satisfies $n>r_A(6+d_\ell)$
there are iid random variables $\eta_1,\eta_2,\dots$  with values in $\{1,2,\dots,[n/r_A]\}$
and independent of the Poissonian $Z$ with parameter $t(1-\rho_A)$ such that 
$$
\mathbb{P}\left(\xi_n^t\in L\right)=\mathbb{P}(W\in L)+\mathcal{E}(A)\hspace{4mm}\forall L\subset\mathbb{N}_0,
$$
where $W=\sum_{i=1}^Z\eta_i$ is the associated compound Poisson random variable.

The error  term $\mathcal{E}(A)$ equals $e^{-\gamma n}+\psi(n)$ up to polynomial terms in $n$
and exponential terms in $\ell,t$, where $\gamma>0$ is such that 
$\mu(B)\le e^{-2\gamma m}\;\forall B\in\mathcal{A}^m\;\forall m\in\mathbb{N}$.
\end{theorem}

\noindent This result completely describes (subject to the condition $n>r_A(6+d_\ell)$)
the distribution at every point, periodic or not. In the case $\ell=1$ and $d_1=1$ at a periodic point $x$ 
of minimal period $m$ we have $r=m$ for $n$ large enough, $\rho=p$ from Theorem~\ref{PolyaAeppli}
 and the random variables $\eta_i$ can be replaced by random variables that are geometrically 
 distributed with parameter $p$. The limiting result is then the same as in Theorem~\ref{PolyaAeppli}
 although the error terms are larger.

%%%%%%%%%%%%%%%%%%%%%%%%%%%%%%%%%%%%%%%%%%
%%%%%%%   EXAMPLES
%%%%%%%%%%%%%%%%%%%%%%%%%%%%%%%%%%%%%%%%%%
\section{Recurrence times}

Recurrence time is the special case when the point whose reentry is observed is the 
same at which the target cylinders are centered. To be precise,
let $\mathcal{A}$ be a generating finite or countably infinite partition of $\Omega$,
then $R_n(x)=\tau_{A_n(x)}(x)$ is the recurrence time and measures the time 
it takes for the first $n$ symbols of a point $x$ to reoccur in $x$. In the symbolic description, when 
    every point $x$ is identified by its trajectory $\vec{x}=(\dots,x_{-1},x_0,x_1,\dots)$ then
  $$
  R_n(x)=\min\{j\ge1: x_jx_{j+1}\cdots x_{j+n-1}=x_0x_{1}\cdots x_{n-1}\}
  $$
   measures  the time it takes until one sees the starting $n$-word again. 
  According to Kac's theorem the value of $\tau_{A_n(x)}$ is on average $1/\mu(A_n(x))$.
  Denote by $h(\mu)$ the measure theoretic entropy of the invariant probability 
  measure $\mu$.
  According to the theorem of Shannon-McMillan-Breiman~\cite{Man} one has $\mu(A_n(x))\sim e^{-nh}$ 
  which would one make expect that
  $R_n(x)\sim e^{nh}$. This indeed is true as was proven by Ornstein and Weiss 
  first for finite alphabets in 1993~\cite{OW1} and in 2002 for countably infinite alphabets~\cite{OW2}:
    
  \begin{theorem} \cite{OW1,OW2} Let $\mu$ be ergodic and $\mathcal{A}$ 
  a finite or countably infinite $\mu$-generating partition, then  almost surely
  $$
  \lim_{n\rightarrow\infty}\frac{\log R_n(x)}n=h(\mu).
  $$
  \end{theorem}
  
   \noindent In the infinite case one must have $\sum_{A\in\mathcal{A}}\mu(A)|\log \mu(A)|<\infty$
   in order to ensure finite entropy.
   
   For some mixing systems this result was strengthened by
Kontoyiannis~\cite{Kon} who prove the almost sure invariance principle.
The requirement is that the invariant measure is $\alpha$-mixing of a sufficient 
rate ($\alpha(\Delta)=\mathcal{O}(\Delta^{-336})$) and satisfy an
$\mathscr{L}^1$-Gibbs condition\footnote{
Denote by $f_n(\vec{x})=-\mathbb{P}(x_0|x_{-1}x_{-2}\cdots x_{-n})$ and let 
$f$ be the pointwise limit of $f_n$ as $n\to \infty$. It is then required that
$\|f-f_n\|_1=\mathcal{O}(n^{-48})$}. 
This strengthened a previous result of  Nobel and Wyner~\cite{NW} who showed that for 
strongly mixing systems (without regularity condition) the exponential growth rate of 
recurrence times equals the metric entropy.

For exponentially $\psi$-mixing Gibbs measures, Collet, Galves and Schmitt~\cite{CGS}
proved the Central Limit Theorem for the recurrence time that is 
$(\log R_n-nh(\mu))/\sigma\sqrt{n}$ converges in distribution to the normal law
 (provided the variance $\sigma^2$, which is given by the Gibbs potential, is positive).
 This required the CLT for 
Shannon-McMillan-Breiman and the fact that entry times are exponentially distributed.

A similar result holds for metric spaces. If $T$ is a map on a metric space $\Omega$ 
with metric $d$, then the $n$th Bowen ball is given by 
$B_{\varepsilon,n}(x)=\{y\in\Omega: d(T^jx,T^jy)<\varepsilon, 0\le j<n\}$. 
With F Yang we have proven that for an ergodic $T$-invariant probability measure $\mu$
one has
$$
\lim_{\varepsilon\to0}\lim_{n\to\infty}\frac1n\log R_{\varepsilon,n}(x)=h(\mu)
$$
almost everywhere, where $R_{\varepsilon,n}(x)=\tau_{B_{\varepsilon,n}(x)}(x)$ is 
the recurrence time to the Bowen ball (the limit in $n$ is $\limsup$ or $\liminf$).

For geometric balls one can define $R_\varepsilon(x)=\tau_{B_\varepsilon(x)}(x)$ for the 
recurrence time to the geometric ball $B_\varepsilon(x)$. Then
$\lim_{\varepsilon\to0}\frac{\log R_\varepsilon(x)}{|\log \varepsilon|}\le d(x)$ if
the limits and the dimension exist and otherwise for $\limsup$ and $\liminf$ on both 
sides. This was shown in~\cite{BS} to be true almost everywhere for invariant measures
of maps on compact manifolds. Equality was proven by Saussol~\cite{Sau06} and 
more generally in~\cite{RS} in the case when correlation functions decay
superpolynomially.

An interesting connection to the R\'enyi entropy function is provided by looking 
at the quantity
$$
Z_n(t)=\sum_{j=0}^{R_n(x)}\mu(A_n(T^jx))^t
$$
and its limiting behaviour as $n$ goes to infinity. For the equilibrium state $\mu$ for H\"older 
continuous potential $\varphi$ on a subshift of finite type its exponential growth rate is 
$$
\lim_{n\to\infty}\frac1n\log Z_n(t)=h(\mu)+P((1+t)\varphi)-(1+t)P(\varphi)
$$
almost everywhere, where $P$ denotes the pressure function. 
For $t=1$ this was shown by  Deschamps, Schmitt, Urbanski and Zdunik~\cite{DSUZ}
by using large deviations which in this case are exponential.
Evidently, for $t=0$ one recovers the statement of Ornstein and Weiss' theorem. 
The general case for $\phi$-mixing measures was dealt with by Ko~\cite{Ko}  who showed that 
$$
\lim_{n\to\infty}\frac1n\log Z_n(t)=h(\mu)-tR(t),
$$
 where
$R(t)=\lim_{n\to\infty}\frac1{nt}\log\sum_{A\in\mathcal{A}}\mu(A)^{1+t}$ is the R\'enyi
entropy function ($R(0)=h(\mu)$).

For the generalised recurrence time we have~\cite{GHKR}
$$
\liminf_{n\to\infty}\frac1n\log\tau_{A_n(z)}(x)\ge h(\mu)
$$
for $\mu\times\mu$ almost all $(z,x)\in\Omega\times\Omega$ for all ergodic invariant
measures $\mu$. Unfortunately a general result as in Ornstein and Weiss' theorem 
cannot hold true as Shields~\cite{Shi} has produced an example of an invariant measure
over a four element shift space for which 
$\limsup_{n\to\infty}\frac1n\log\tau_{A_n(z)}(x)=\infty$ almost surely with respect to the product
measure.
For Gibbs measures on subshifts of finite type it Chazottes and Ugalde~\cite{CU} have proven
that $\frac1n\log\tau_{A_n(z)}(x)\to h(\mu)$ almost surely in $\mu\times\mu$.
More generally for  $\phi$-mixing measures it was shown in~\cite{GHKR}  that
under the assumption
that the limiting entry times function $F_x(t)=\lim_{n\to\infty}F_x^n(t)$ exists almost surely
and decays to $0$ as $t$ goes to infinity, then 
$$
\frac1n\log\tau_{A_n(z)}(x)\longrightarrow h(\mu)
$$
 in measure.

%%%%%%%%%%%%%%%%%%%%%%%%%%%%%%%%%%%%%%%%
%%%%%%%%%%%  REFERENCES
%%%%%%%%%%%%%%%%%%%%%%%%%%%%%%%%%%%%%%%%

\end{document}